\documentclass{amsart}

\usepackage{a4wide,vaucanson-g}

\newtheorem{lemma}{Lemma}[section]
\newtheorem{thm}[lemma]{Theorem}
\newtheorem{definition}[lemma]{Definition}
\newtheorem{proposition}[lemma]{Proposition}
\newtheorem{corollary}[lemma]{Corollary}

\def\R{\mathcal R}
\def\El{\mathcal E_\ell}
\def\Er{\mathcal E_r}

\title{Sequences with constant number of return words}

\author{L\!'ubom\'ira Balkov\'a}
\author{Edita Pelantov\'a}
\address{Doppler Institute for Mathematical Physics and Applied
Mathematics, and Department of Mathematics, FNSPE, Czech Technical
University, Trojanova 13, 120~00 Praha~2, Czech Republic}
\email{l.balkova@centrum.cz, Pelantova@km1.fjfi.cvut.cz}

\author{Wolfgang Steiner}
\address{LIAFA, CNRS, Universit\'e Paris Diderot -- Paris 7, Case 7014, 
75205 Paris Cedex 13, France}
\email{steiner@liafa.jussieu.fr}
\date\today

\begin{document}
\begin{abstract}
An infinite word has the property $R_m$ if every factor has exactly
$m$ return words. 
Vuillon showed that $R_2$ characterizes Sturmian words. 
We prove that a word satisfies $R_m$ if its complexity function is
$(m-1)n+1$ and if it contains no weak bispecial factor.
These conditions are necessary for $m=3$, whereas for $m=4$ the 
complexity function need not be $3n+1$. 
New examples of words satisfying $R_m$ are given by words related to 
digital expansions in real bases.
\end{abstract}

\maketitle

\section{Introduction}
Recently, return words have been intensively studied in (symbolic) 
dynamical systems, combinatorics on words and number theory. 
Roughly speaking, for a~given factor $w$ of an infinite word $u$, 
a~return word of $w$ is a~word between two successive occurrences of 
the factor~$w$. 
This can be seen as a symbolic version of the first return map in a 
dynamical system.
This notion was introduced by Durand~\cite{Durand} to give a~nice 
characterization of primitive substitutive sequences. 
A~slightly different notion of return words was used by Ferenczi, 
Mauduit and Nogueira~\cite{FeMaNo}. 

Sturmian words are aperiodic words over a~biliteral alphabet with the
lowest possible factor complexity; they were defined by Morse and
Hedlund~\cite{MoHe}.
Using return words, Vuillon~\cite{Vuillon} found a~new equivalent
definition of Sturmian words.
He showed that an infinite word $u$ over a~biliteral alphabet is
Sturmian if and only if any factor of $u$ has exactly two return words.
A short proof of this fact is given in Section~\ref{R2R3}.

A natural generalization of Sturmian words to $m$-letter alphabets is 
constituted by infinite words with every factor having exactly $m$ 
return words. 
This property is called $R_m$.
It covers other generalizations of Sturmian words:
Justin and Vuillon~\cite{JusVui} proved that Arnoux-Rauzy words of
order $m$ satisfy $R_m$, Vuillon~\cite{VuIET} proved this property for 
words coding regular $m$-interval exchange transformations.

The factor complexity, i.e., the number of different factors of length 
$n$, of the two classes of words with property $R_m$ in the preceding 
paragraph is $(m-1)n+1$ for all $n\ge0$. 
Vuillon~\cite{VuIET} observed that this condition is not sufficient to 
describe words satisfying $R_m$, $m\ge3$: the fixed point of a certain
recoding of the Chacon substitution, which has complexity $2n+1$ by 
Ferenczi~\cite{Fer}, has factors with more than $3$ return words.

A deeper inspection of the two classes of words with property $R_m$
shows that not only the first difference of complexity is constant, but 
also that the bilateral order of every factor (see 
Cassaigne~\cite{Cass} and Section~\ref{sufficient}) is zero. 
We show that this condition is indeed sufficient to have the property
$R_m$, and provide a less known class of words satisfying this
condition.
If a word satisfies $R_3$, then we can show that no factor is weak
bispecial, i.e., no factor has negative bilateral order.
Therefore the words with $R_3$ are characterized by complexity $2n+1$
and the absence of weak bispecial factors. 

In Section~\ref{cex}, we provide a word satisfying $R_4$ with an even 
number of factors of every positive length (containing infinitely many 
weak bispecial factors). 
Therefore words satisfying $R_m$ do not necessarily have complexity 
$(m-1)n+1$, and it is an open question whether there exists a nice 
characterization of words satisfying $R_m$ for $m\ge4$.

We conclude the article by exhibiting a large class of purely 
substitutive words satisfying $R_m$. 
Every word in this class codes the sequence of distances between 
consecutive $\beta$-integers for some real number $\beta>1$.

In this article we focus only on the number of return words 
corresponding to a given factor of an infinite word. 
We do not study the ordering of return words in the infinite word, 
i.e., we do not study derivated sequences (see~\cite{Durand} for the 
precise definition). 
Let us just mention here that a derivated sequence of a word with 
property $R_m$ is again a word satisfying $R_m$.
A~description of derivated sequences of Sturmian words can be found 
in~\cite{ArBr}.


\section{Basic definitions} \label{Preliminaries}

An {\em alphabet} $\mathcal A$ is a~finite set of symbols called
{\em letters}. A~(possibly empty) concatenation of letters is a~{\em
word}. The set $\mathcal A^{*}$ of all finite words provided with
the operation of concatenation is a~free monoid. The {\em length} of
a~word $w$ is denoted by $|w|$. A~finite word $w$ is called a~{\em
factor} (or {\em subword}) of the (finite or right infinite) word
$u$ if there exist a~finite word $v$ and a~word $v'$ such that
$u=vwv'$. The word $w$ is a~{\em prefix} of $u$ if $v$ is the empty
word. Analogously, $w$ is a~{\em suffix} of $u$ if $v'$ is the empty
word. A~concatenation of $k$ words $w$ will be denoted by $w^k$.

The {\em language} $\mathcal{L}(u)$ is the set of all factors of the
word $u$, and $\mathcal{L}_n(u)$ is the set of all factors of $u$ of
length $n$.
Let $w$ be a~factor of an infinite word $u$ and let $a,b\in\mathcal A$.
If $aw$ is a~factor of $u$, then we call $a$ a~{\em left extension} of
$w$.
Analogously, we call $b$ a~{\em right extension} of~$w$ if
$wb\in\mathcal L(u)$.
We will denote by $\El(w)$ the set of all left extensions of $w$, and
by $\Er(w)$ the set of right extensions.
A~factor $w$ is {\em left special} if $\#\El(w)\ge2$, {\em right
special} if $\#\Er(w)\ge2$ and {\em bispecial} if $w$ is both left
special and right special.

Let $w$ be a~factor of an infinite word $u=u_0u_1\cdots$ (with
$u_j\in\mathcal A$), $|w|=\ell$.
An integer $j$ is called an {\em occurrence} of $w$ in $u$ if
$u_ju_{j+1}\cdots u_{j+\ell-1}=w$.
Let $j,k$, $j<k$, be successive occurrences of $w$.
Then $u_ju_{j+1}\cdots u_{k-1}$ is a~{\em return word} of $w$.
The set of all return words of $w$ is denoted by $\R(w)$,
$$
\R(w)=\{u_ju_{j+1}\dots u_{k-1}\mid j,k \mbox{ being successive
occurrences of } w \mbox{ in }u\}.
$$
If $v$ is a return word of $w$, then the word $vw$ is called
{\em complete return word}.

An infinite word is {\em recurrent} if any of its factors occurs
infinitely often or, equivalently, if any of its factors occurs at
least twice.
It is {\em uniformly recurrent} if, for any $n\in\mathbb N$, every 
sufficiently long factor contains all factors of length $n$.
It is not difficult to see that a recurrent word on a finite alphabet
is uniformly recurrent if and only if the set of return words of any
factor is finite.

The variability of local configurations in $u$ is expressed by the
{\em factor complexity function} (or simply {\em complexity})
$ C(n)=\#\mathcal L_n(u)$. It is well known that a~word $u$ is
aperiodic if and only if $ C(n)\ge n+1$ for all $n \in \mathbb N$.
Infinite aperiodic words with the minimal complexity $ C(n)=n+1$ for
all $n\in \mathbb N$ are called {\em Sturmian words}. These words
have been studied extensively, and several equivalent definitions of
Sturmian words can be found in Berstel~\cite{Berstel}. 


\section{Simple facts for return words}\label{HandyRules}

\subsection{Restriction to bispecial factors}\label{sectbisp}
If a~factor $w$ is not right special, i.e., if it has a~unique right
extension $b\in\mathcal A$, then the sets of occurrences of $w$ and
$wb$ coincide, and
$$
\R(w)=\R(wb).
$$
If a~factor $w$ has a~unique left extension $a\in\mathcal A$, then
$j\ge 1$ is an occurrence of $w$ in the infinite word $u$ if and only
if $j-1$ is an occurrence of $bw$. This statement does not hold for
$j=0$.
Nevertheless, if $u$ is a recurrent infinite word, then the set of
return words of $w$ stays the same no matter whether we include the
return word corresponding to the prefix $w$ of $u$ or not.
Consequently, we have
$$
\R(aw)=a\R(w)a^{-1}=\{ava^{-1}\mid v\in\R(w)\},
$$
where $ava^{-1}$ means that the word $v$ is prolonged to the left by
the letter $a$ and it is shortened from the right by erasing the
letter $a$ (which is always the suffix of $v$ for $v\in\R(w)$).

For an aperiodic uniformly recurrent infinite word $u$, each factor
$w$ can be extended to the left and to the right to a~bispecial factor.
To describe the cardinality and the structure of $\R(w)$ for
arbitrary $w$, it suffices therefore to consider bispecial factors $w$.

\subsection{Tree of return words}\label{tree}
It is convenient to consider a tree (or trie) constructed in the 
following way:
Label the root with a factor $w$, and attach $\#\Er(w)$ children,
with labels $wb$, $b\in \Er(w)$.
Repeat this recursively with every node labeled by $v$, except if $w$
is a suffix of $v$.
If $u$ is uniformly recurrent, then this algorithm stops, and it is
easy to see that the labels of the leaves of this tree are exactly the
complete return words of $w$.
Therefore we have
\begin{equation}\label{leaves}
\#\R(w)=\#\{\text{leaves}\}=1+\sum_{\text{non-leaves }v}(\#\Er(v)-1).
\end{equation}
In particular, if $w$ is the unique right special factor of its length,
then $\#\R(w)=\#\Er(w)$.

\begin{figure}
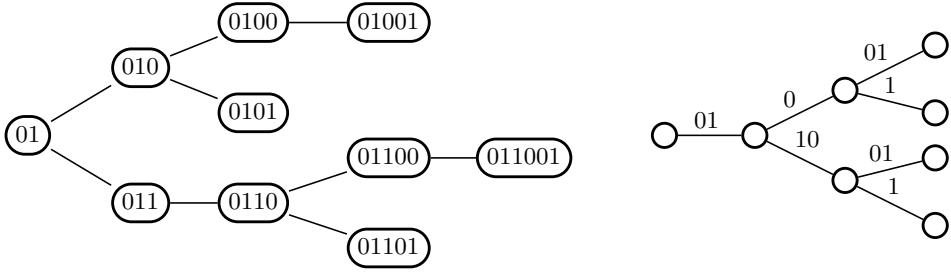

{\small\VCDraw{\begin{VCPicture}{(0,-3)(20,3)}
\StateVar[01]{(0,0)}{a}
\StateVar[010]{(2.5,1.5)}{b}
\StateVar[011]{(2.5,-1.5)}{c}
\StateVar[0100]{(5,2.5)}{d}
\StateVar[0101]{(5,.5)}{e}
\StateVar[0110]{(5,-1.5)}{f}
\StateVar[01001]{(8,2.5)}{g}
\StateVar[01100]{(8,-.5)}{h}
\StateVar[01101]{(8,-2.5)}{i}
\StateVar[011001]{(11,-.5)}{j}

\SetEdgeArrowStyle{-}
\EdgeL ab{}
\EdgeL ac{}
\EdgeL bd{}
\EdgeL be{}
\EdgeL cf{}
\EdgeL dg{}
\EdgeL fh{}
\EdgeL fi{}
\EdgeL hj{}

\SmallState
\State{(14,0)}{k}
\State{(16,0)}{l}
\State{(18,1)}{m}
\State{(18,-1)}{n}
\State{(20,2)}{o}
\State{(20,.5)}{p}
\State{(20,-.5)}{q}
\State{(20,-2)}{r}

\EdgeL kl{01}
\EdgeL lm{0}
\EdgeL ln{10}
\EdgeL mo{01}
\EdgeL mp{1}
\EdgeL nq{01}
\EdgeL nr{1}
\end{VCPicture}}}
\caption{The tree of return words of $01$ in the Thue-Morse sequence
and its trie representation.}
\end{figure}

A similar construction can be done with left extensions, yielding
similar formulae.
Since we can restrict our attention to bispecial factors $w$ by
Section~\ref{sectbisp}, we obtain the following proposition.

\begin{proposition}\label{onespecial}
Let $u$ be a recurrent word and $m\in\mathbb N$.
Suppose that for every $n\in\mathbb N$ at least one of the following
conditions is satisfied:
\begin{itemize}
\item
There is a unique left special factor $w\in\mathcal L_n(u)$, and
$\#\El(w)=m$.
\item
There is a unique right special factor $w\in\mathcal L_n(u)$, and
$\#\Er(w)=m$.
\end{itemize}
Then $u$ satisfies property $R_m$, i.e., every factor has exactly $m$
return words.
\end{proposition}

Recall that Arnoux-Rauzy words of order $m$ are defined as uniformly
recurrent infinite words which have for every $n\in\mathbb N$ exactly
one right special factor $w$ of length $n$ with $\#\Er(w)=m$ and
exactly one left special factor $w$ of length $n$ with $\#\El(w)=m$.
They are also called strict episturmian words.
It is easy to see that Sturmian words are recurrent, and we obtain the
following corollary to Proposition~\ref{onespecial}.

\begin{corollary}\label{ArnouxRauzy}
Arnoux-Rauzy words of order $m$ satisfy $R_m$, in particular Sturmian 
words satisfy $R_2$.
\end{corollary}

\section{Sufficient conditions for property $R_m$}\label{sufficient}

This section is devoted to sufficient conditions for a word $u$ having 
the property $R_m$, but we mention first two evident necessary 
conditions.

The alphabet $\mathcal A$ of $u$ must have $m$ letters since the 
occurrences of the empty word are all integers $n\ge0$, and its return 
words are therefore all letters $u_n$. 
Furthermore, $u$ must be uniformly recurrent since every factor has a 
return word and only finitely many of them.

An important role in our further considerations is played by weak
bispecial factors.

\begin{definition}\label{weak}
A factor $w$ of a recurrent word is {\em weak bispecial} if $B(w)<0$,
where 
$$
B(w)=
\#\{awb\in\mathcal L(u)\mid a,b\in\mathcal A\}-\#\El(w)-\#\Er(w)+1
$$
is the {\em bilateral order} of $w$.
\end{definition}

Since
$\#\{awb\in\mathcal L(u)\mid a,b\in\mathcal A\}= 
\sum_{a\in \El(w)}\#\Er(aw)=\sum_{b\in \Er(w)}\#\El(wb)$,
the inequality $B(w)<0$ is equivalent to
$$
\sum_{a\in \El(w)}(\#\Er(aw)-1)<\#\Er(w)-1
$$
and to
$$
\sum_{b\in \Er(w)}(\#\El(wb)-1)<\#\El(w)-1.
$$

The bilateral order was defined by Cassaigne~\cite{Cass} in order to 
calculate the second complexity difference.
If we set $\Delta C(n)=C(n+1)-C(n)$, then we have 
$$
\Delta C(n) = \sum_{w \in \mathcal{L}_n(u)} \bigl(\#\El(w)-1\bigr)
= \sum_{w \in \mathcal{L}_n(u)} \bigl(\#\Er(w)-1\bigr)
$$
and therefore
\begin{multline*}
\Delta C(n+1)-\Delta C(n)=\sum_{w\in\mathcal L_n(u)}\sum_{a\in \El(w)}
(\#\Er(aw)-1)-\sum_{w\in\mathcal L_n(u)}(\#\Er(w)-1) \\
=\sum_{w\in\mathcal L_n(u)}\big(\#\{awb\in\mathcal L(u)\mid a,b\in
\mathcal A\}-\#\El(w)-\#\Er(w)+1\big)=\sum_{w\in\mathcal L_n(u)}B(w).
\end{multline*}

If $B(w)=0$ for all factors $w$, then the first complexity difference 
is constant.
If no factor is weak bispecial, then $\Delta C(n)$ is non-decreasing.
Since $\Delta C(0)=\#\mathcal A-1$ and $\#A=m$, we obtain the following 
lemma.

\begin{lemma}\label{l0}
If $u$ satisfies $R_m$ and no factor is weak bispecial, then 
$\Delta C(n)\ge m-1$ for all $n\ge0$.
\end{lemma}


The number of return words can be bounded by the following lemmas.

\begin{lemma}\label{l1}
If $u$ is a uniformly recurrent word with no weak bispecial factor,
then
$$
\#\R(w)\ge1+\Delta C(|w|)
$$
for every factor $w\in\mathcal L(u)$.
\end{lemma}

\noindent{\em Proof.} Let $w\in\mathcal L(u)$ and denote by
$v_1,v_2,\ldots,v_r$ the right special factors of length $|w|$.
Since no factor is weak bispecial and $u$ is uniformly recurrent,
every $v_j$ can be extended to the left without decreasing the total
amount of ``right branching'' until $w$ is reached. 
More precisely, we have (mutually different) right special factors
$v_j^{(1)},v_j^{(2)},\ldots,v_j^{(s_j)}$ with suffix $v_j$, prefix
$w$ and no other occurrence of $w$ such that
$\#\Er(v_j)-1\le\sum_{i=1}^{s_j}(\#\Er(v_j^{(i)})-1)$. Since all
$v_j^{(i)}$ are nodes in the tree of return words and $v_j^{(i)}\ne
v_{j'}^{(i')}$ if $(j,i)\ne(j',i')$, we can use (\ref{leaves}) and
obtain
$$
\qquad\#\R(w)\ge1+\sum_{j=1}^r\sum_{i=1}^{s_j}(\#\Er(v_j^{(i)})-1)\ge
1+\sum_{j=1}^r(\#\Er(v_j)-1)=1+\Delta C(|w|).\qquad\qed
$$

\begin{lemma}\label{l2}
If $u$ has no weak bispecial factor and $\Delta C(n)<m$ for all 
$n\ge0$, then
$$
\#\R(w)\le m
$$
for every factor $w\in\mathcal L(u)$.
\end{lemma}

\noindent{\em Proof.} Let $v_1,v_2,\ldots,v_r$ denote the right
special factors which are labels of non-leave nodes in the tree of
return words of $w$, and $n=\max_{1\le j\le r}|v_j|$. Since no 
bispecial factor is weak, every $v_j$ can be extended to the left to
factors of length $n$ without decreasing the total amount of ``right
branching''. More precisely, we have (mutually different) right
special factors $v_j^{(1)},v_j^{(2)},\ldots,v_j^{(s_j)}$ of length
$n$ with suffix $v_j$ such that
$\#\Er(v_j)-1\le\sum_{i=1}^{s_j}(\#\Er(v_j^{(i)})-1)$. Since $w$
occurs in $v_j$ only as prefix, no $v_j$ can be a proper suffix of
$v_{j'}$. Hence we have $v_j^{(i)}\ne v_{j'}^{(i')}$ if
$(j,i)\ne(j',i')$ and
$$
\#\R(w)=1+\sum_{j=1}^r\big(\#\Er(v_j)-1\big)\le 1+\sum_{j=1}^r
\sum_{i=1}^{s_j}\big(\#\Er(v_j^{(i)})-1\big)\le1+\Delta C(n)\le m.\qed
$$

For words with no weak bispecial factors, these three lemmas give a 
very simple characterization of the property $R_m$.

\begin{thm}\label{th1}
If $u$ is a uniformly recurrent word with no weak bispecial factor, 
then it satisfies $R_m$ if and only if $ C(n)=(m-1)n+1$ for all $n\ge0$.
\end{thm}

\section{Properties $R_2$ and $R_3$}\label{R2R3}

For $m=2$ and $m=3$, we can completely characterize the words with
property $R_m$.

\begin{definition}
Let $v$ be a return word of $w\in\mathcal L(u)$. We say that return 
word $v$ starts with $b$ if $wb$ is a prefix of the complete return 
word $vw$ and that it ends with $a$ if $aw$ is a suffix of $vw$.
\end{definition}

A right special factor $w$ is called {\em maximal right special} if $w$ 
is not a proper suffix of any right special factor, i.e.,
$\sum_{a\in \El(w)}(\#\Er(aw)-1)=0$. 
Any maximal right special factor is therefore weak bispecial.

\begin{lemma}\label{lemma1}
If $w\in\mathcal L(u)$ is a maximal right special factor such that
for any $b\in \Er(w)$ there exists a unique $v\in\R(w)$ starting
with $b$, then $u$ is eventually periodic.
\end{lemma}

\begin{proof}
Denote the return words of $w$ by $v_1,v_2,\ldots,v_r$, where,
w.l.o.g., $v_j$ starts with $b_j$, ends with $a_j$ and $b_{j+1}$ is
the only letter in $\Er(a_jw)$ for $1\le j<r$. Then $b_1$ is the
only letter in $\Er(a_rw)$ and $u=p(v_1v_2\cdots v_r)^\infty$ for
some prefix $p$.
\end{proof}

\begin{corollary}\label{cormax}
If $u$ satisfies $R_2$, then it has no maximal right special factor.
\end{corollary}

\begin{proof}
Assume that $w$ is a maximal right special factor. Then the two
return words of $w$ have different starting letters, hence $u$ is
eventually periodic by Lemma~\ref{lemma1} and $\#\R(wa)=1$.
\end{proof}

On a binary alphabet, the notions ``weak bispecial'' and ``maximal
right special'' coincide. 
Therefore Corollaries~\ref{ArnouxRauzy}, \ref{cormax} and 
Lemma~\ref{l1} provide a short proof of the following theorem.

\begin{thm}[Vuillon~\cite{Vuillon}]
An infinite word $u$ satisfies $R_2$ if and only if it is Sturmian.
\end{thm}

For words with property $R_3$, we need the following lemma.

\begin{lemma}\label{lemma2}
Let $w$ be a weak bispecial factor with a unique $a\in \El(w)$ such
that more than one return word of $w$ starts with a letter in
$\Er(aw)$, then $\#\R(aw)<\#\R(w)$.
\end{lemma}

\begin{proof}
Any return word of $aw$ has the form $av_1v_2\cdots v_ra^{-1}$ for
some $r\ge1$ and $v_j\in\R(w)$, $1\le j\le r$. If $v_1$ ends with
$a$, then $r=1$. If $v_1$ ends with $a'\ne a$, then the assumption
of the lemma implies that there is a unique return word of $w$
starting with a letter in $\Er(a'w)$ (and $\#\Er(a'w)=1$).
Therefore $v_2$ and inductively the sequence of words
$v_2,\ldots,v_r$ are completely determined by the choice of $v_1$.
This implies that $\#\R(aw)$ equals the number of return words of
$w$ starting with a letter in $\#\Er(aw)$. Since $w$ is weak 
bispecial, we have $\#\Er(aw)<\#\Er(w)$ and thus $\#\R(aw)<\#\R(w)$.
\end{proof}

\noindent{\it Remark.}\
There are two cases for Lemma~\ref{lemma2}: 
Either $aw$ is right special or there is more than one return word of 
$w$ starting with the unique right extension of $aw$.

\begin{corollary}\label{nomaximal}
If $u$ satisfies $R_3$, then it has no weak bispecial factor.
\end{corollary}

\begin{proof}
Assume that $w$ is a weak bispecial factor.

If $\#\Er(w)=3$, then every return word of $w$ starts with a
different letter in $\Er(w)$. Since at most for one $a\in \El(w)$, 
the factor $aw$ is right special, we obtain a contradiction to $R_3$ by 
Lemma~\ref{lemma1} or \ref{lemma2}.

If $\#\Er(w)=2$, then $\Er(aw)=\{b\}$ and $\Er(a'w)=\{b'\}$.
Since, w.l.o.g., two return words of $w$ start with $b$ and one
starts with $b'$, we obtain a contradiction to $R_3$ by
Lemma~\ref{lemma2}.
\end{proof}

By combining  Corollary~\ref{nomaximal} and Theorem~\ref{th1}, we
obtain the following theorem.

\begin{thm}
A uniformly recurrent word $u$ satisfies $R_3$ if and only if
$ C(n)=2n+1$ for all $n\ge0$ and $u$ has no weak  bispecial factor.
\end{thm}

\noindent{\it Remarks.}
\begin{itemize}
\item
The theorem remains true if ``weak bispecial'' is replaced by 
``maximal right special'': 
If $\Delta C(n)=2$ for all $n\ge0$, then every factor $w$ with 
$\#\Er(w)=3$ is the unique right special factor of its length, and it 
cannot be weak bispecial. 
If $\#\Er(w)=2$, then the two notions coincide.
\item
By symmetry, ``weak bispecial'' can be replaced by ``maximal left 
special''.
\item
The condition on weak bispecial factors cannot be omitted.
Ferenczi~\cite{Fer} showed that the fixed point $\sigma^\infty(1)$
of the substitution given by 
$\sigma:1\mapsto12,2\mapsto312,3\mapsto3312$, a recoding of the
Chacon substitution, has complexity $2n+1$ and it contains weak 
bispecial factors.
\end{itemize}

\section{Property $R_4$}\label{R_4}

\subsection{A word with complexity $\ne3n+1$}\label{cex}
The following proposition shows that $C(n)$ need not be $(m-1)n+1$ for 
all $n\ge0$ if $u$ satisfies $R_m$.

\begin{proposition}
Define the substitution $\sigma$ by
\begin{align*}
\sigma: 1&\mapsto 13231 \\
2&\mapsto 13231424131 \\
3&\mapsto 42324131424 \\
4&\mapsto 42324
\end{align*}
Then the fixed point $\sigma^\infty(1)$ satisfies $R_4$.
\end{proposition}

\begin{proof}
By Section~\ref{sectbisp}, it is sufficient to consider bispecial 
factors of $u=\sigma^\infty(1)$.
The factors of length $2$ are 
$\mathcal L_2(u)=\{13,14,23,24,31,32,41,42\}$.
For the bispecial factors $1,2,23,2413$, the return words are easily
determined:
\begin{align*}
\R(1) & = \{13,1323,1424,142324\} \\
\R(2) & = \{23,2314,2413,241314\} \\
\R(23) & = \{2314,2314241314,232413,232413142413\} \\
\R(2413) & = \{241314,24131423,24132314,2413231423\} 
\end{align*}

The language of $u$ is closed under the morphism $\varphi$ defined by 
$\varphi:1\leftrightarrow4,2\leftrightarrow3$, since 
$\sigma\varphi(w)=\varphi\sigma(w)$ for all factors $w$.
Therefore we have $\R(\varphi(w))=\varphi(\R(w))$.

The only factors of the form $a1b$, $a,b\in\mathcal A$ are $314$ and 
$413$, hence $1$ is a weak bispecial factor, and $1,4$ are the only 
bispecial factors with prefix or suffix $1$ or $4$.
Similarly, $23$ and $32$ are weak bispecial factors and no other 
bispecial factor has prefix or suffix $23$ or $32$.

The return words of the weak bispecial factor $2413142$ are factors of 
$\sigma(v)$, with a factor $v$ of length $|v|\ge2$ having prefix $2$ or $3$, 
suffix $2$ or $3$ and no other occurrence of $2$ and $3$.
Since the only possibilites for $v$ are $23,2413,32,3142$, we obtain
$$
\R(2413142)= 
\{24131423,241314232413231423,24131424132314,241314241323142324132314\}.
$$

All remaining bispecial factors $w$ have prefix $24132$ or $31423$ and 
suffix $23142$ or $32413$, and therefore a decomposition 
$w=t\,\sigma(v)\,t'$ with $t\in\{24,31\}$, $t'\in\{1323142,4232413\}$ 
and a unique bispecial factor $v$.
If $v$ is empty, then we have w.l.o.g. $w=241323142$ and
$$
\R(w)=\{2413231423,2413231423241314,2413231424131423,
24132314241314232413242\}.
$$
If $v$ is not empty, then the uniqueness of $v$ implies that the set of 
complete return words of $w$ is $t\,\sigma(\R(v)v)\,t'$. 
Since $v$ is shorter than $w$, we obtain inductively that all bispecial
factors have exactly $4$ return words.
\end{proof}

\subsection{Weak bispecial factors}

The preceding example shows that weak bispecial factors cannot be 
excluded in words $u$ satisfying $R_4$. 
Nevertheless, we can show that the existence of a weak bispecial 
factor imposes strong restrictions on the structure of the word $u$.

\begin{lemma}
Let $w$ be a weak bispecial factor of a word $u$ satisfying $R_4$.
Then there exist factors $w_1,w_2\in\mathcal Aw\cup w\mathcal A$ and
$v_1,v_2,v_3,v_4$ such that
\begin{equation}\label{ret1234}
\R(w_1)=\{v_1v_3,v_1v_4,v_2v_3,v_2v_4\}\mbox{ and }
\R(w_2)=\{v_3v_1,v_3v_2,v_4v_1,v_4v_2\}.
\end{equation}
\end{lemma}

\begin{proof}
Let $w$ be a weak bispecial factor. Since the problem is symmetric,
assume, w.l.o.g., $\#\Er(w)\ge\#\El(w)$. We have three different
situations:
\begin{itemize}
\item
$\#\Er(w)=2$: We have $\Er(aw)=1$ for all $a\in \El(w)$. For both
$b\in \Er(w)$, there must be more than one return word of $w$
starting with $b$ by Lemma~\ref{lemma2}. Therefore, there exists two
return words of $w$ starting with each $b$. Let $w_1=wb_1$,
$w_2=wb_2$.
\item
$\#\Er(w)=3$: There exists a unique $b\in \Er(w)$ such that two
return words of $w$ start with~$b$. By Lemma~\ref{lemma2}, there
exists therefore some $a\in \El(w)$ such that $\#\Er(aw)\ge2$.
Since $w$ is weak bispecial with $\#\Er(w)=3$, we have $\#\Er(aw)=2$
and $\Er(a'w)=1$ for all other $a'\in \El(w)$.
Let $w_1=aw$, $w_2=wb$.
\item
$\#\Er(w)=4$: For every $b\in \Er(w)$, there is a unique return
word of $w$ starting with~$b$. By Lemma~\ref{lemma2}, we have
$a_1,a_2\in \El(w)$ with $\#\Er(a_iw)=2$. Let $w_1=a_1w$,
$w_2=a_2w$.
\end{itemize}

Consider ``complete return words of the set $\{w_1,w_2\}$'': words
which have either $w_1$ or $w_2$ as prefix, either $w_1$ or $w_2$ as
suffix, and no other occurrence of $w_1$ and $w_2$. By the
definitions of $w_1$ and $w_2$, there are exactly two such words
$v_1w_{i_1},v_2w_{i_2}$ with prefix $w_1$ and two words 
$v_3w_{i_3},v_4w_{i_4}$ with prefix $w_2$.

If $i_1=i_2=2$ and $i_3=i_4=1$, then $R_4$ implies that (\ref{ret1234}) 
holds.

If $i_1=i_2=1$, then $w_1$ has only the two return words $v_1,v_2$. 
If $i_2=i_3=i_4=1$, then the return words of $w_1$ are 
$v_1v_3,v_1v_4,v_2$. 
Similarly, $i_3=i_4=2$ and $i_1=i_2=i_3=2$ are not possible.

The only remaining case is $i_1=i_4=1$, $i_2=i_3=2$. 
Then the return words of $w_1$ are $v_1$ and $v_2v_3^{r_i}v_4$,
$i\in\{1,2,3\}$, $0\le r_1<r_2<r_3$. The return words of $w_2$ are
$v_3$ and $v_4v_1^{s_i}v_2$, $i\in\{1,2,3\}$, $0\le s_1<s_2<s_3$.

The return words of $v_2w_2$ are therefore of the form
$v_2v_3^{r_i}v_4v_1^{s_j}$. Let $S_1$ be the set of these $4$ pairs
$(r_i,s_j)$. Similarly, let $S_2$ be the set of the $4$ pairs
$(s_i,r_j)$ such that $v_4v_1^{s_i}v_2v_3^{r_j}$ is a return word of
$v_4w_1$.

We show that there must be some $i\in\{1,2,3\}$ such that
$(r_i,s_2)\in S_1$ and $(r_i,s_3)\in S_1$, by considering the return
words of $v_1^{s_2}w_1$ and of $v_1^{s_2}v_2w_2$. The return words
of $v_1^{s_2}v_2w_2$ are of the form
$v_1^{s_2}v_2tv_3^{r_i}v_4v_1^{s_j-s_2}$ with
$t\in(v_3^*v_4v_1^{s_1}v_2)^*$, $i\in\{1,2,3\}$ and $j\in\{2,3\}$.
For these $t$ and $r_i$, $v_1^{s_2}v_2tv_3^{r_i}v_4$ is a return
word of $v_1^{s_2}w_1$. If there was no $r_i$ with $(r_i,s_2)\in
S_1$ and $(r_i,s_3)\in S_1$, then these words would provide $4$
different return words of $v_1^{s_2}w_1$, wich contradicts $R_4$
since $v_1$ is another return word.

Similarly, we must have some $i\in\{1,2,3\}$ such that $(s_i,r_2)\in
S_2$ and $(s_i,r_3)\in S_2$. By considering the return words of
$v_1^{s_2}w_1$ and $v_4v_1^{s_2}w_1$, we obtain as well the
existence of some $i\in\{1,2,3\}$ such that $(r_2,s_i)\in S_1$ and
$(r_3,s_i)\in S_1$. Finally, we must also have some $i\in\{1,2,3\}$
such that $(s_2,r_i)\in S_2$ and $(s_3,r_i)\in S_2$.

Putting everything together, we have two possibilities for $S_1$.
Either it contains $(r_1,s_1)$ and no other $(r_i,s_j)$ with $i=1$
or $j=1$, or $S_1=\{(r_1,s_2),(r_1,s_3),(r_2,s_1),(r_3,s_1)\}$.
Analogously, $S_2$ contains $(s_1,r_1)$ and no other $(s_i,r_j)$
with $i=1$ or $j=1$, or
$S_2=\{(s_1,r_2),(s_1,r_3),(s_2,r_1),(s_3,r_1)\}$.

If $(r_1,s_1)\in S_1$ and $(s_1,r_1)\in S_2$, then
$v_2v_3^{r_1}v_4w_1$ has only one return word,
$v_2v_3^{r_1}v_4v_1^{s_1}$. If $(r_1,s_1)\not\in S_1$ and
$(s_1,r_1)\not\in S_2$, then $v_2v_3^{r_1}v_4w_1$ has only two
return words, $v_2v_3^{r_1}v_4v_1^{s_2}$ and
$v_2v_3^{r_1}v_4v_1^{s_3}$. If $(r_1,s_1)\in S_1$ and
$(s_1,r_1)\not\in S_2$, then the return words of
$v_2v_3^{r_1}v_4w_1$ are of the form
$v_2v_3^{r_1}v_4v_1^{s_1}v_2v_3^{r_i}v_4v_1^{s_j}$ with
$(r_i,s_j)\in S_1\setminus\{(r_1,s_1)\}$, thus there are only three
words. Similarly, $v_4v_1^{s_1}v_2w_2$ has only three return words
if $(r_1,s_1)\not\in S_1$ and $(s_1,r_1)\in S_2$.

This shows that $i_1=i_4=1$, $i_2=i_3=2$ is impossible, and
the lemma is proved.
\end{proof}

\section{Words associated with $\beta$-integers}\label{beta}

In this section, we describe a new class of infinite words with
the property $R_m$. The language of these words is not necessarily
closed under reversal.

Consider the fixed point $u=\sigma^\infty(0)$ of a primitive
substitution of the form
\begin{equation}\label{subst}
\begin{array}{rcl}\sigma:\qquad 0 & \mapsto & 0^{t_1}1 \\
1 & \mapsto & 0^{t_2}2 \\ & \vdots \\
m-2 & \mapsto & 0^{t_{m-1}}(m-1) \\ m-1 & \mapsto & 0^{t_m}\end{array}
\end{equation}
for some integers $m\ge2$, $t_1,t_m\ge1$ and $t_2,\ldots,t_{m-1}\ge0$.
The incidence matrix of $\sigma$ is a companion matrix of the
polynomial $x^m-t_1x^{m-1}-\cdots-t_m$.
Let $\beta>1$ be the dominant root of this polynomial (the
Perron-Frobenius eigenvalue of the matrix).
If
$$
t_j\cdots t_m\prec t_1\cdots t_m\quad\mbox{for all }j\in\{2,\ldots,m\},
$$
where $\preceq$ denotes the lexicographic ordering, then $\sigma$ is a
{\em $\beta$-substitution} and $\beta$ is a simple Parry number.
It is easy to see that $u$ codes in this case the sequence of
distances between consecutive nonnegative $\beta$-integers
$$
\mathbb Z_\beta^+=\Big\{\sum_{j=0}^J x_j\beta^j\mid J\ge0,\,
x_j\in\mathbb Z,\, x_j\ge0,\, x_j\cdots x_0\prec t_1\cdots t_m
\text{ for all }j,\,0\le j\le J\Big\},
$$
and a letter $k$ corresponds to the distance
$t_{k+1}/\beta+\cdots+t_m/\beta^{m-k}$.
($0$~corresponds to distance~$1$.)

\medskip\noindent{\em Remark.}
The most prominent example of a $\beta$-substitution is the Fibonacci
substitution ($m=2$, $t_1=t_2=1$), where $\beta$ is the golden mean.

\medskip
It is not difficult to show that all prefixes of $u$ are left special
factors, with all $m$ letters being left extensions (see
e.g. Frougny, Mas\'akov\'a and Pelantov\'a~\cite{FrMaPe}).
For every factor $w$, the tree of return words constructed by the left
extensions (see Section~\ref{tree}) contains therefore a node with $m$
children, the shortest prefix of $u$ having $w$ as suffix.
The word $u$ is uniformly recurrent since all fixed points of primitive
substitutions have this property (Queff\'elec~\cite{Queffelec}).
Therefore every factor $w$ has at least $m$ return words.
If there exists a left special factor which is not a prefix of $u$,
then this factor has more than $m$ return words.
By Proposition~\ref{onespecial}, we obtain the following proposition.

\begin{proposition}
If $u=\sigma^\infty(0)$ for some substitution $\sigma$ of the form
(\ref{subst}), then it satisfies $R_m$ if and only if 
$ C(n)=(m-1)n+1$ for all $n\ge0$.
\end{proposition}

Bernat, Mas\'akov\'a and Pelantov\'a~\cite{BeMaPe} characterized the
fixed points of $\beta$-substitutions satisfying $\Delta C(n)=m-1$ for 
all $n\ge0$.
The techniques of their proof can also be used to construct non-prefix 
left special factors if $\sigma$ is a substitution of the form 
(\ref{subst}) which is not a $\beta$-substitution, and their conditions
can be reformulated as in the following corollary.

\begin{corollary}\label{affineRm}
If $u=\sigma^\infty(0)$ for some substitution $\sigma$ of the form
(\ref{subst}), then it has the property $R_m$ if and only if
\begin{itemize}
\item $t_m=1$ and
\item $t_j\cdots t_{m-1}t_1\cdots t_{j-1}\preceq t_1\cdots t_{m-1}$
for all $j\in\{2,\ldots,m-1\}$.
\end{itemize}
\end{corollary}

Note that the language of $u$ is closed under reversal if and only if
$t_1=t_2=\cdots=t_{m-1}$.
In this case, $u$ is an Arnoux-Rauzy word of order $m$.


\section*{Acknowledgements}

The authors acknowledge financial support by the Czech Science
Foundation GA \v{C}R 201/05/0169 and by the grant LC06002 of the
Ministry of Education, Youth and Sports of the Czech Republic.


\end{document}